\input amstex
\documentstyle{amsppt}
\input bull-ppt

%\define\index{\operatorname{index}}

%\redefine\qed{\hfill \vrule \height5truept \width3truept 
%\depth0truept}
%\redefine\qed{\hfill $\blacksquare$}
\define\ddt{\frac{d}{dt}}
%\define\({\left(}
%\define\){\right)}
\NoBlackBoxes

\topmatter
\cvol{26}
\cvolyear{1992}
\cmonth{Jan}
\cyear{1992}
\cvolno{1}
\cpgs{125-130}
\title Nonunique tangent maps at isolated\\
singularities
of harmonic maps
\endtitle
\shorttitle{Nonunique Tangent Maps}
\author Brian White \endauthor
\affil Stanford University \endaffil
\address Mathematics Department, Stanford University, 
Stanford,
California 94305\endaddress
\thanks The author was partially funded by NSF grants 
DMS85-53231 (PYI),
DMS87-03537, and DMS9012-718\endthanks
\abstract Shoen and Uhlenbeck showed that ``tangent 
maps'' can be defined at
singular points of energy minimizing maps. Unfortunately 
these are not unique,
even for generic boundary conditions. Examples are 
discussed which have
isolated singularities with a continuum of distinct 
tangent maps.
\endabstract
\date April 3, 1991\enddate
\endtopmatter
\document

Let $\Omega$ be a bounded domain in $R^m$ (or more 
generally a compact
riemannian manifold with boundary) and let $N$ be a 
compact riemannian
manifold.  By the Nash embedding theorem, $N$ can be 
regarded as a
submanifold of some euclidean space.
The energy of a map $f\:\Omega\to N$ is defined to be
$$
  E(f\,) = \int_{\Omega} |Df|^2.
$$
(Here $f$ is allowed to be any measurable map from 
$\Omega$ to $R^d$ such that
$f(x)\in N$ for almost every $x$ and such that the 
distributional first
derivative of $f$ is square integrable.)
The map $f$ is said to be energy minimizing if its energy 
is less than
or equal to the energy of each other map having the same 
boundary values.
It is fairly easy to prove that if $g\:\Omega \to N$ has 
finite energy, then
there is an energy minimizing map $f\:\Omega \to N$ with 
the same
boundary values as $g$.
In \cite{SU}, Schoen and Uhlenbeck proved that if $f$ is 
energy minimizing,
then $f$ is smooth except on a set $K\subset \Omega$ of 
Hausdorff
dimension at most $m-3$.

Suppose $f$ is energy minimizing and that $x\in \Omega$ 
is a singularity of
$f$.  Schoen and Uhlenbeck also proved that for every 
sequence $r_i$ of
positive numbers converging to zero, a subsequence
of the maps
$$
        y \mapsto f(x+r_iy)    \tag 1
$$
converges weakly to a map $f_\infty\:R^m\to \Omega$ that 
is constant
on rays through the origin.  Such a map is called a {\it 
tangent map\/}
to $f$ at $x$.   Intuitively, $f_\infty$ is the result of 
looking at $f$
near $x$ through a microscope with infinite 
magnification.  The map $f_\infty$
is simpler than $f$ because it is constant on rays, but 
one would like to
think that it provides a good picture of $f$ near $x$. 
  Note that $f_\infty$
would not give a very good picture of $f$ if there were 
more than one
tangent map at $x$; that is, if a different subsequence 
of the maps \thetag1
could converge to another limit map.   Whether or not 
such pathological
behavior is possible has been perhaps the most basic open 
question about
singularities of energy minimizing maps.

There have been some positive results (ruling out 
pathological behavior).
First, Leon Simon \cite{S} showed that if $N$ is analytic 
and if $f_\infty$ has
an isolated discontinuity, then $f_\infty$ is unique 
(i.e., it is the
only tangent map at $x$).  Second, Gulliver and White 
\cite{GW} showed that if
$m=3$ and $\dim(N)=2$
(the lowest dimensions in which singularities are 
possible),
then $f_\infty$ is unique whether or not $N$ is analytic.

This paper is an announcement of the first example of 
nonuniqueness:

\proclaim{Theorem 1}
There exists a $C^\infty$ $5$-manifold $N$ and a nonempty 
open
set $U$ of smooth maps $\phi:\partial B^4\to N$ such that
\roster
\item each $\phi\in U$ bounds one or more energy 
minimizing maps from $B^4$
to $N$, and
\item if $f\:B^4\to N$ is an energy minimizing map
with $f\mid\partial B^4\in U$,
then $f$ has an isolated singularity $x$ and a continuum 
of tangent maps at
$x$. Each of the tangent maps is regular except at $0$.
\endroster
\endproclaim

The proof is too long to give here; see \cite{W4}.
However, we can prove a simpler but
nonetheless interesting result that has the same flavor. 
 Let $f\:\Omega\to N$
be a finite energy map that is smooth except at a finite 
set of discontinuities
$\{p_i:i=1,\dots,k\}$.  We say that $f$ is {\it harmonic} 
if it satisfies
the Euler-Lagrange partial differential equations for the 
energy functional.
Such an $f$ is a critical point for energy, but it need 
not be a minimum.
Nonetheless, the existence of tangent maps and the 
uniqueness results of
Simon and Gulliver-White in fact hold for such harmonic 
maps.
Thus it is interesting to note:

\proclaim{Theorem 2}
There is a $C^\infty$ $4$-manifold $N$ and a harmonic map 
$f\:B^3\to N$
such that $f$ has an isolated singularity at $0$ and a 
continuum of
distinct tangent maps at $0$.
\endproclaim

\demo{Proof} Let $N$ be the product $S^1\times R\times 
S^2$ with
the metric
$$
      dx^2 + dy^2 + (2-V(x,y))dz^2
$$
$(V$ will be specified later).
Note that this defines a complete metric on $N$ provided 
$V$ is everywhere
less than $2$.
Of course $N$ is not compact, but the image
of the harmonic map we construct will be contained in a 
compact subset of
$N$, so we could easily modify $N$ to make it compact.

Note that the orthogonal map $O(3)$ acts on $B^3$ and on 
$N$ (on $N$ by
$\rho\:(x,y,z)\mapsto (x,y,\rho z)$).  We simplify the 
harmonic map equations
by looking for solutions that are 
$O(3)$-equivariant.  It is not hard to see
that every equivariant map is of the form:
$$
   p \mapsto  (v^1(|p|),v^2(|p|),\pm p/|p|) \in S^1\times 
R\times S^2.\tag 2
$$
Here $v\:(0,1]\to S^1\times R$.  It is convenient to 
introduce a change
of variable.  Let $t=\log r$ and $u(t)=v(e^t)$, so
$u\:(-\infty,0]\to S^1\times R$.
Then the energy of the map \thetag{2} is
$$
      \int_{t=-\infty}^0 (|\dot u|^2 + (2-V(u)))e^{t}\,dt.
$$
(If the domain were $k$-dimensional, then $e^t$ would be 
$e^{(k-2)t}$.)
The associated Euler-Lagrange equation is
$$
     \ddot u + \dot u + \nabla V(u) = 0.  \tag 3
$$

Thus equivariant harmonic maps are equivalent to 
solutions of the
ordinary differential equation \thetag{3}.  This equation 
has a
physical interpretation: it is the equation of motion of 
a unit mass
moving in $S^1\times R$ subject to a potential $V$ and a 
viscous
force.  Thus the physical energy $E(u,t)=\tfrac12|\dot 
u|^2+V(u)$ is
monotonically decreasing.  To see this mathematically, 
multiply
\thetag{3} by $\dot u$:
$$
        \ddt \left(\frac12 \dot u^2 + V(u)\right) = -\dot 
u^2.  \tag 4
$$

Now we choose $V$ to be
$$
   V(x,y) = -\exp\!\left(-\frac1{y^2}\right)\sin\!\left(x+
\frac1y\right).
$$

\proclaim{Lemma}
$(1)$ Let $u$ be a solution of \thetag{3} with 
$E(u,\cdot)$ constant.
Then $u(t)\equiv p$ for some $p\in S^1\times [0]$.

$(2)$ Let $u$ be a solution of \thetag{3} such that the
physical energy $E(u,0)$ at time $0$ is negative.  Then 
the solution
exists for all $t\in [0,\infty)$ and becomes unbounded as 
$t\to\infty$.
\endproclaim
\demo{Proof}
If $E(u,\cdot)$ is constant, then $\dot u\equiv 0$, since 
otherwise
the particle would be dissipating physical energy to 
viscosity (see
\thetag{4}).  Thus $u$ is a constant $p$, so \thetag{3} 
implies that
$\nabla V(p)=0$.  But $\nabla V(x,y)=0$ if and only if 
$y=0$.
This proves \therosteritem1.

Now suppose that $u$ is a solution of \thetag{3} with 
$E(u,0)<0$.
By elementary ODE theory, the solution exists for all 
positive
times unless the particle moves infinitely far in a 
finite time.
But $\frac12 |\dot u(t)|^2 + V(u(t)) \le E(u,0) < 0$,
so $\frac12 |\dot u(t)|^2 < -V(u(t)) \le \sup(-V)=1$. 
 Thus the solution
exists for $t\in [0,\infty)$.

Suppose $u(t)$ remains in a bounded region of $S^1\times 
R$.
Then the set of pairs $(u(n),\dot u(n))$ is bounded, so a 
subsequence
$(u(n_i),\dot u(n_i))$ converges.  It follows (from the 
smooth dependence of
ODE solutions on initial conditions) that
the solutions $u_i(t)=u(n_i+t)$ converge smoothly to a 
solution
$v(t)$ of \thetag{3}.  Now
$$
 E(v,t) = \lim_{i\to \infty} E(u,n_i+t)
        = \lim_{t\to \infty}E(u,t) \le E(u,0) < 0 \tag 5
$$
(where $\lim_{t\to \infty}E(u,t)$ exists because 
$E(u,\cdot)$ is monotonic).
Thus $E(v,\cdot)$ is constant, so by \therosteritem1 of 
the lemma
$v(t)\equiv p \in S^1\times [0]$.  But then 
$E(v,t)=V(p)=0$, contradicting
\thetag{5}.  This proves \therosteritem2. \qed\enddemo

Now let $u_n\:[0,\infty)\to S^1\times R$ be the solution 
to \thetag{3} with
initial position
$u_n(0) = (\tfrac12\pi, \tfrac1{2\pi n})
$
and initial velocity $\dot u_n(0)=0$.
Note that the initial physical energy is negative:
$$
 \frac12|\dot u_n(0)|^2 + V(u_n(0)) = 0 +V(\frac12\pi, 
\frac1{2\pi n})<0.
$$
Thus by the lemma, there is a first time $t_n>0$ at which
$u_n(t_n)\in S^1\times \{-1,1\}$.
If $u_n(t)$ were ever in $S^1\times [0]$, then 
$E(u_n,t)\ge V(u_n(t))=0$,
which is impossible. Thus $u_n(t_n) \in S^1\times [1]$ and
$u_n(t)\in S^1\times (0,1)$ for $t\in (0,t_n)$.

Note that $u_n(0)\to (\tfrac12\pi,0)$ and $\dot 
u_n(0)\equiv 0$, so
the $u_n$ converge to a solution $w$ with 
$w(0)=(\tfrac12\pi,0)$ and
$\dot w(0)=0$.  By uniqueness of solutions to ODEs,
$w(t)\equiv (\tfrac12\pi,0)$.  Thus
$\lim_{n\to \infty}u_n(t)=(\tfrac12\pi,0)$, so $t_n\to 
\infty$
since $u_n(t_n)=(x_n,1)$.

Now as in the proof of the lemma, there is a sequence 
$n(i)$
such that the solutions $v_i(t) = u_{n(i)}(t_{n(i)}+t)$ 
converge smoothly
to a solution $v$ on $(-\infty,\infty)$.  Of course
$$
\gather
       v(0)\in S^1\times [1], \\
       \text{$v(t) \in S^1\times [0,1]$\quad for $\quad 
t<0$},\\
\intertext{and}
E(v,t)\le 0\quad\text{for all } t.
\endgather
$$
In fact $E(v,t)$ must be strictly negative for every $t$. 
 For since
$E(v,t)$ is a nonpositive and nonincreasing function of 
$t$, if it were $0$
for some $t=a$, then it would be $0$ for each $t\le a$. 
 But then by the
lemma, $v(t)\equiv p \in S^1\times [0]$ for all $t\le a$. 
 By unique
continuation for ODE, $v(t)\equiv p \in S^1\times [0]$ 
for all $t$.
But $v(0)\in S^1\times [1]$.  This proves that $E(v,t)$ 
is strictly
negative.

Now I claim that $v$ defines a harmonic map with a 
continuum of tangent
maps at the origin.  That is, I claim that $v(t)$ has a 
continuum of
subsequential limits as $t\to -\infty$.

As in the proof of the lemma, every sequence of $t$'s 
tending to $-\infty$ has
a subsequence $\tau_i$ such that the solutions 
$w_i(t)=v(t+\tau_i)$ converge to
a solution $w(t)$.
Of course
$$
  E(w,t)=\lim_{i\to\infty}E(v,t+\tau_i) = 
\lim_{t\to-\infty}E(v,t) \le 0
$$
(where $\lim_{t\to -\infty}E(v,t)$ exists because 
$E(v,\cdot)$ is monotonic).
Thus $E(w,\cdot)$ is constant, so by the lemma 
$w(t)\equiv p$, where $p\in
S^1\times [0]$.

What we have shown is $\lim_{t\to -\infty}v^2(t) = 0$, 
where $v^2(t)$ is the
second component of $v(t)=(v^1(t),v^2(t)) \in S^1\times R$.

Now the set $Z=\{p\in S^1\times R: V(p)=0\}$
consists of $S^1\times [0]$ together with a collection of 
curves that
wind around the cylinder infinitely many times as they 
approach
$S^1\times[0]$.  Since $V(v(t)) \le E(v,t) < 0$, $v(t)$ 
is never in $Z$.
Thus $v(t)$ must also wind around the cylinder infinitely 
many times as
$t\to-\infty$.  This proves Theorem 2.

(To make this last argument more formal, note from the 
definition of $V$
that for each $x\in S^1$ and each $\varepsilon>0$, the set
$Z\cup ([x]\times (-\varepsilon,
\varepsilon))$ divides $S^1\times R$ into infinitely
many connected components, the closure of each of which 
is disjoint from
$S^1\times [0]$.  Since $v(t)$ approaches $S^1\times [0]$ 
as
$t\to -\infty$, the particle must cross the set
$Z\cup ([x]\times (-\varepsilon,\varepsilon))$.  Since it 
never crosses $Z$, it
must cross $[x]\times (-\varepsilon,
\varepsilon)$.  As this holds for every $x$ and
$\varepsilon$,
 each $(x,0)\in S^1\times [0]$ is a subsequential limit 
of $v(t)$.)
\qed\enddemo

\heading Remarks \endheading
%\subheading{Remarks}
Exactly the same construction provides examples of 
harmonic maps from
$B^m$ to $N=S^1\times R \times S^{m-1}$ (metrized as 
above) with
a continuum of tangent maps at an isolated singularity.
The only difference is that the viscosity (i.e., the 
coefficient in front
of $\dot u$ in \thetag{3}) is $m-2$ instead of $1$.

In all those examples, the dimension of the target 
manifold is one more than
the dimension of the domain.
But we can also prove that there is a harmonic map $f$ 
from $B^4$ to the
$4$-manifold $N$ of Theorem 2 such that $f$ has a 
continuum of tangent maps
at an isolated singularity.  The proof is the same as the 
proof of Theorem 2,
except that we consider maps of the form
$$
    f\: p \to (f^1(|p|), f^2(|p|), h(p/|p|)),
$$
where $h:S^3\to S^2$ is the Hopf fibration.

\heading Open questions \endheading
%\subheading{Open questions}
\subheading{1} Must tangent maps be unique if the target 
manifold $N$ is
$2$-dimensional?  The answer is ``yes'' if the domain is 
$3$-dimensional
\cite{GW}.
\subheading{2} Must tangent maps be unique for generic 
metrics on the target
manifold $N$?
\subheading{3} If $T$ is a minimal variety in a 
riemannian manifold $N$, then
at each singular point $x \in T$ there are one or more 
tangent cones
(i.e.,  subsequential limits of images of $T$ under 
dilations about $x$).
Can there be more than one?  See \cite{AA; T1,2; W1--3},
and Simon \cite{S} for results in special cases.
Simon \cite{S} proved that if a tangent cone has 
multiplicity one and
has an isolated singularity, then it is unique.  Unlike 
his analogous
result for harmonic maps, this does not require that the 
metric on $N$
be analytic.

The construction in this paper does not seem to have any 
analogue in the case
of minimal varieties.

\enddocument